\RequirePackage{fix-cm}
\documentclass[aps,showkeys,showpacs]{revtex4}
\usepackage{amsmath,amsfonts,amssymb,amscd, euscript}
\usepackage[dvips]{graphicx}
\begin{document}

\title{Effective algorithm of analysis of integrability via the Ziglin's method.
}


\author{Vladimir Salnikov}


\affiliation{              Laboratoire de Math\'ematiques de l'INSA de Rouen, \\ 
Avenue de l'Universit\'e 76801 Saint-\'Etienne-du-Rouvray Cedex, France \\
              Tel.: +33 2 32 95 65 52          
}
              \email{vladimir.salnikov@insa-rouen.fr}

\begin{abstract}
In this paper we continue the description of the possibilities to use numerical 
simulations for mathematically rigorous computer assisted analysis of integrability
of dynamical systems.
We sketch some of the algebraic methods of studying the integrability and 
present a constructive algorithm issued from the Ziglin's approach. 
We provide some examples of successful applications of the constructed algorithm 
to physical systems. 
\end{abstract}
\maketitle

\section{Introduction}
In this paper we describe the continuation of the research 
aimed at application of numerical methods to analysis of integrability 
of dynamical systems. 
In the first paper of the series (\cite{int-num}) we have described 
the way of revealing the obstructions to real integrability in the \emph{Liouville--Arnold 
sense} (\cite{AKN}) by analyzing the topology of the phase space of the system. 
We have also discussed the possibilities of extending the approach 
to parametrized families of dynamical systems with the goal of 
searching for the regions of possible integrability.

The current paper is devoted to analysis of algebraic properties of integrable systems. 
We mention some recent results related to the study of complexified systems and systems of
\emph{variational equations} in particular. We pay special attention to the results of S.L.~Ziglin (\cite{Ziglin}) on the \emph{meromorphic integrability}, that is for the system with $d$ degrees of freedom 
we are interested in
the existence of $d$ independent first integrals in involution, that are meromorphic 
functions of the phase space coordinates. 
We point out the major difficulties of application of these results to the study of integrability 
and propose an effective algorithm of a computer assisted construction significantly
extending the range of their applicability.
We give some examples of application of the method to dynamical systems 
having physical origin, and also mention some possible purely mathematical 
outcome. 

\section{Algebraic obstructions to integrability}
Since for a given dynamical system there is no general approach for studying 
the existence of the sufficient number of arbitrary first integrals, a 
natural idea to restrict the class of first integrals comes out. 
A rather detailed review of the appropriate methods can be found 
in \cite{kozlov_sym}, some more recent methods are also well explained in 
 \cite{audin}. Here we will only sketch some of them, that are important to understand 
 our motivations and the results presented in this paper. 

One can probably say, that the first algebraic method for analysis of integrability 
is the approach of S.~Kovalevskaya developed in \cite{yoshida}
by H.~Yoshida for studying the polynomial first integrals of 
dynamical systems. Let us note that this restriction is rather reasonable 
since for mechanical systems with polynomial hamiltonian functions the natural integrals 
(energy, angular momentum etc.) are of this class.
H.~Yoshida applies this method to study the Euler equations as well as some hamiltonian systems 
having physical origin.

\subsection{Variational equations, monodromy group.} \label{monodromie}
Some more recent methods of analysis of integrability deal with the complexified systems
of differential equations. The key observation is that rather often for classical 
completely integrable systems the first integrals can be continued to the 
complex domain of the canonical variables depending on complex time remaining in 
the class of holomorphic or meromorphic functions. 
And branching of the solutions of the hamiltonian equations can create an obstruction 
to the existence of such first integrals. This idea permits in particular to 
find the relations between the parameters of a system necessary for complete integrability, 
that is select more regular systems from the family of a priori  similar ones. 

Following partially \cite{audin} let us formalize the above statement.
The following approach actually dates back to A.~Lyapunov who proposed to study the 
system of variational equations, it was developed by S.~Ziglin in \cite{Ziglin}
who revealed the relation of the structure of the monodromy group 
to integrability. 

Consider in $\mathbb C^n$ (or any $n$-dimensional complex manifold) a system
\begin{equation}
\label{2syst1}
  \dot x = v(x),
\end{equation}
where the dot denotes the derivative with respect to the complex time.
Denote $x_0(t)$ -- a particular solution of (\ref{2syst1}) and consider a 
small perturbation $x = x_0 + \xi$ of it. 
Note that here and in what follows the variables like $x, v$ will be 
vector-valued (in $\mathbb C^n$, in a complex manifold or in the 
appropriate tangent bundle), 
and those like $\xi$ can be identified to a respective tangent vector field
in the whole phase space.
Plugging $x$ into 
(\ref{2syst1}) one obtains
$$
\dot \xi = A\xi + \dots, \qquad A = \frac{\partial v}{\partial x}(x_0(t)).
$$
Then the linearized system 
\begin{equation}
\label{syst_lin}
 \dot \xi = A(x_0(t))\xi
\end{equation}
is called the system of {\it variational equations} along the particular solution $x_0(\cdot)$.
We write $A(x_0(t))$ instead of just $A(t)$, to stress the fact that the matrix of this linear system 
depends on time only implicitly via its dependence on the particular solution.
If the system (\ref{syst_lin}) has multivalued solutions then the system 
(\ref{2syst1}) also does. 
We will be however interested in a more subtle property of these systems, namely 
branching of the solutions of the system of variational equations (\ref{syst_lin})  
when moving along the Riemann surface of the particular solution  $x_0$ of (\ref{2syst1}).

Thus we will consider a system (\ref{syst_lin}) of $n$ linear equations, 
where the entries of the matrix 
$A$ are holomorphic functions defined in a connected neighbourhood
of the Riemann surface $x_0$.
Locally for any given value of $\xi(t_0) = \xi_0$ there exists a unique holomorphic 
continuation of the solution (\ref{syst_lin}). One can continue it along any path in $x_0$
but the result in the generic case need not be single-valued. 
Let $\gamma$ be a loop (closed path) starting from $x_0(t_0)$
parametrized by some path on the complex plane. Constructing the continuation 
of the solution  $\xi(t)$ of (\ref{syst_lin}) defined in the neighbourhood of $x_0(t_0)$ along $\gamma$ we obtain another solution  $\xi_*(t)$ of (\ref{syst_lin}) defined in the same neighbourhood.
Since the system (\ref{syst_lin}) is linear, 
there exists a fundamental system of its solutions, such a continuation corresponds 
to a change of basis in the linear space spanned by it and thus can be encoded by a 
complex $n \times n$ matrix $T_\gamma$.
Branching of the solutions of (\ref{syst_lin}) corresponds to 
$T_\gamma \neq id$.

The set of all matrices  ${\cal M} = \{T_\gamma\}$, corresponding to all the loops 
$\gamma$ in $x_0$, forms a group called the \emph{monodromy group}
of the linear system (\ref{syst_lin}) along the solution $x_0$.
The matrices certainly depend on the choice of the base point on $x_0$
but all the groups ${\cal M}$ are isomorphic.
It is important to note that two homotopic loops produce the same 
monodromy matrices, therefore it is sufficient to consider only one 
representative $\gamma$ from each homotopy class.

\subsection{Results of S.~Ziglin \label{sympl} } 
Let  $f(x)$  be a first integral of  (\ref{2syst1}).
One can write down the Taylor series for $f(x_0 + \xi)$ 
$$
  f(x_0 + \xi) = \sum \limits_{m \ge 0} F_m | _{x_0}(\xi, t),
$$
where $F_m$ is the homogeneous form of degree $m$, single-valued on the Riemann surface $x_0$.
It is clear that a non-zero homogeneous form 
$F_m$ of the lowest degree ($m \ge 1$) is the first integral of the system of variational equations 
(\ref{syst_lin}), and thus invariant under the action of the monodromy group.
$$
  F_m(T\xi,t_0) = F_m(\xi,t_0), \quad T \in {\cal M}.
$$
If $f(x)$ is meromorphic, it can be represented as the ratio of two holomorphic functions: 
$\frac{P(z)}{Q(z)}$. Then the analogous invariant is the ratio of the lowest degree 
forms corresponding to $P$ and $Q$: $\frac{P_m(z)}{Q_k(z)}$.
This is a strong condition on invariant functions, but even more important, 
the existence of such invariants imposes serious restrictions on the structure of 
the monodromy group.

Without going much into details let us just mention that for what follows it is important 
to exclude the cases when all the matrices of the monodromy group have trivial eigenvalues. 
This happens when the system possesses so-called symmetry fields; if this is the case
one needs to perform the reduction of the system to normal variational equations. 
Let us also note that for the case of hamiltonian systems, which is of particular 
interest for us in the context of analysis of integrability, the elements of the 
monodromy group define symplectic affine transformations. In this case 
the eigenvalues of the monodromy matrices split into couples of mutually inverse 
complex numbers. 
The transformation is called \emph{non-resonant} if its eigenvalues 
 ($\lambda_1, \lambda_1^{-1}, \dots, \lambda_{p}, \lambda_{p}^{-1}$)
satisfy the relation
 $
   \lambda_1^{k_1} \cdot \dots \cdot \lambda_{p}^{k_{p}} = 1
 $
if and only if all the $k_i$'s vanish.

Having defined the notion of the monodromy group we can state the following theorems:

{\bf Theorem 1. (Ziglin's lemma)}
\textit{
Let the monodromy group of the curve  $x_0$ contain a non-resonant transformation 
$g$.
Then the number of meromorphic first integrals of the hamiltonian equations 
in the connected neighbourhood of $x_0$ functionally independent from $H$
is not bigger than the number of rational invariants of the monodromy group.
}

{\bf Theorem 2. (Ziglin)}
\textit{
Let the monodromy group ${\cal M}$ of the curve  $x_0$ contain a non-resonant transformation 
$g$.
If the hamiltonian equations admit ($d-1$) 
meromorphic first integrals in the connected neighbourhood
of $x_0$ functionally independent from $H$ then 
any 
other transformation $g' \in {\cal M}$ preserves the fixed point of $g$
and maps its eigendirections to eigendirections. If moreover no set of 
eigenvalues of $g'$ forms a regular polygon on the complex plane 
(with the number of vertices $\ge 2$), $g'$ preserves the eigendirections of $g$
(i.e. commutes with $g$).
}

\subsection{Differential Galois theory}
Despite some important examples, the analytical computation of the monodromy 
group is a complicated task, solved only for some particular systems. 
In this context it is worth noting that there is an alternative approach 
proposed by J.~Morales-Ruiz and J.-P.~Ramis (\cite{morales,morales-ramis})
following the scheme of S.~Ziglin. 
It establishes the relation between the properties of the differential Galois group
 and integrability. Not defining this group here (for the definition see e.g. \cite{shafar}) 
 let us only state the main result. 

{\bf Theorem 3. (Morales-Ramis)} 
\textit{Let the hamiltonian system be completely integrable, then the 
identity component of the differential Galois group of the system of normal variational 
equations along any particular solution is abelian.
}

When restricted to the case of complete integrability, i.e. the existence of a sufficient number 
of independent first integrals \emph{in involution}, the theorem of Ziglin is the consequence of
the one of Morales--Ramis. The latter is considered to be stronger (\cite{audin}) also in the sense 
that it permits to show non-integrability in some cases when the Ziglin's approach 
did not give any answer. 
The explanation of this fact is that for a fixed particular solution of (\ref{2syst1})
the differential Galois group contains the monodromy group, and therefore 
has more sources of non-commutativity. 
In this paper we however restrict our attention to the Ziglin's approach, leaving the 
results of Morales--Ramis to further studies.

\section{Effective algorithm of application of the Ziglin's method.}

We can formulate the standard way to analyze the integrability which is 
essentially based on the results reviewed in the previous section.
The major steps are the following:
\begin{itemize}
 \item[1.] For a given complexified system, of differential equations (\ref{2syst1})
 construct explicitly a solution $x_0$, which as a function of complex time is viewed
 as a Riemann surface. 
 \item[2.] Write down the system of variational equations along $x_0$. Perform a reduction 
 of this system to the system of normal variational equations using known first integrals or symmetry fields. 
 \item[3.] Localize the singularities of the particular solution $x_0$.
 \item[4.] Construct the monodromy matrices, obtained by going along the loops around the 
 singularities obtained in 3. -- they generate the monodromy group.
 \item[5.] Make a conclusion about the presence of the obstruction to integrability, based 
 on the commutation relations between the matrices obtained in 4.
\end{itemize}

Let us comment on some subtleties of application of this method. 
The method is obviously not intended to prove integrability, 
and there are 
examples (\cite{morales-ramis,audin}) when the monodromy group is trivial, but the system is still 
non-integrable. 
But it can be used to single out the relations between the parameters 
of the system, when integrability is possible. It is usually done, when the 
variational equations reduce to some ``classical'' well-studied systems.
In principal, one can argue for integrability, when the monodromy group is 
commutative for \emph{any} particular solution, but this is not worth the efforts, 
since knowing explicitly all the solutions of the initial system (\ref{2syst1})
one can perform a more detailed qualitative study of it.
So the method can be developed mainly to search for the obstructions
to integrability.

We can easily see two major difficulties of application of this method. 
First, one needs to construct an explicit particular solution of 
a system of differential equations. And second, this solution should be on 
one hand simple enough, so that the the monodromy group could be computed, and 
on the other hand, it should be non-trivial, so that the computed group 
could contain sufficient number of sources of non-commutativity.
We employ numerical methods to overcome these difficulties, more precisely 
the idea is to compute numerically the generators of the monodromy group
along a particular solution which is also obtained numerically. 
Let us mention that the idea to use numerical methods in the context 
of the Ziglin's method is certainly very natural, however there are very few 
successful implementations of it. We can mention a couple of 
articles (\cite{yoshida2,maciejewski}) on the subject where the authors 
compute the monodromy group numerically based on an explicitly known particular solution:
a straight-line one and using Jacobi elliptic functions respectively. 
Our approach permits to extend the range of applicability of the 
Ziglin's result since we don't have any a priory restriction 
on the form of the particular solution. 

For a given value of complex time $t$ the value of the solution $x(t)$ can be computed numerically 
without any difficulty. One needs just to remember, that the obtained particular solution 
$x(t)$ should be considered not as a function of a point $t \in \mathbb{C}$ 
in a complex plane, but as a function of (the homotopy class of) a path going to $t$
from the initial point $t_0$.
We can always write down the right hand sides of the variational equations  
$\dot \xi = A(x(t)) \xi$ with the matrix $A$ depending on an arbitrary particular solution 
$x$. The difficulty is that if we don't have the analytical expression 
of $x(t)$ we can not plug it into the variational equations, but with the above remark 
it does not make much sense in the general case. 
That is why we use a natural approach of solving the variational equations in parallel 
to the initial system along any path that interests us. 

So, using the observations about the structure of the monodromy group from the previous section 
(\ref{monodromie}) let us solve numerically the following system
\begin{eqnarray} \label{bigsyst}
  \dot x = v(x) \nonumber \\
  \dot \Xi = A(x) \Xi,
\end{eqnarray}
where the first line is identical to the initial complexified system of $n$ differential 
equations (\ref{2syst1}) ,  
the second one is the matrix equation with $A$ being the matrix of the system of variational 
equations depending explicitly on $x$, and $\Xi$ -- unknown $n \times n$ matrix.
For initial data we take $x(t_0) = x_0$ -- an arbitrary point in the phase space and 
 $\Xi(t_0) = id_n$ -- unit matrix.
Going around the loop on the solution $x$ 
we obtain in $\Xi$ precisely the monodromy matrix corresponding to this loop.
Thus going around all the singularities of the 
solution $x$ we can construct the whole monodromy group.

It is important to understand that these loops are paths on the Riemann surface of the solution 
and not just on the complex plane. This is the major difficulty 
that one faces when applying the procedure: given the form of the equations 
(\ref{bigsyst}) we are forced to integrate them against the complex time, 
but the object of interest for us is the behavior of the solution $x$.
More precisely, when we go along a path in the complex plane and neither  $x$ nor $\Xi$
have non-zero variation, this path does not produce a non-trivial monodromy generator. 
If $\Xi$ has non-zero variation, one needs to check that  the corresponding 
values of $x(t)$ have returned to the initial value. 
Only in this case $\Xi$ is the generator of the monodromy group.
If $x(t)$ has not returned to its initial value one needs to continue 
following the path. This difficulty is related to the fact that 
in general the solution $x$ is not a single-valued function of the complex time
$t$, i.e. going around the loop in $\mathbb{C}$ does not always 
produce a closed loop on $x$ and the parametrization used in the 
original work of S.~Ziglin should be considered as 
the parametrization of the loop on $x$.
It is also clear that if $x$ does not return to its initial value after a finite 
number of going around a loop in $\mathbb{C}$ (that is the topology of the singularity 
is logarithmic), such a loop does not correspond to any matrix in the monodromy group.
But such an infinite branching by itself can be an obstruction to, say, 
the existence of first integrals analytic in the phase space coordinates,
at least if the corresponding energy level is compact in the phase space. 
If all the branching points of $x$ are of finite order (such a system 
satisfies the \emph{generalized Painlev\'e property}), then the outcome of the procedure 
is the set of monodromy group generators. Having that, to prove non-integrability 
it is enough to find a couple of non-commuting matrices among them.

There are several issues also worth being commented on. 
First, there is a difficulty which is more technical than conceptual, it is related to the 
localization of the singularities of the solution $x$ (step 3 in the method above). 
Certainly it is impossible to go around all the loops in $\mathbb{C}$, so we have to restrict the analysis to some compact domain
and go around the points of some finite grid. 
That is we do not try to construct the whole monodromy group, but only 
its subgroup, which is however usually enough to reveal the obstruction to integrability
(cf. the example of \cite{maciejewski}).
Second, we have not discussed here the issue of resonant transformations. 
But, given the symplectic nature of the monodromy transformations (cf. section \ref{sympl}) this problem 
arises only from relatively large size of the system (\ref{2syst1}) in question:
for example with two couples of complex eigenvalues the non-resonant
transformation corresponds to non-collinearity of two vectors which is an open condition so it can 
be guaranteed by a numerical test.
In this setting let us also note that one does not have to perform the reduction (in step 2
of the method) of the system of variational equations (\ref{syst_lin}), it is only necessary to know that it can be performed.
It is also important to mention that since the final action (step 5) is also the verification of an open condition 
of the commutator non-vanishing, it is perfectly correct from the point of view of the accuracy of numerical integration.

Summing up,  let us formulate the \textbf{effective algorithm} of analysis of integrability 
via the Ziglin's method.
\begin{enumerate}
 \item[i.] Write down (analytically) the system of equations (\ref{bigsyst}) not fixing the 
 particular solution.
 \item[ii.] Choose a bounded domain in $\mathbb{C}$ and a finite grid of points in it with a 
 distinguished point $t_0$.
 \item[iii.] For each point choose a loop going around only it and starting at $t_0$. 
   Integrate numerically the system (\ref{bigsyst}) along this loop taking 
 $\Xi(t_0) = id$ as the initial conditions. Three cases are possible:
     \begin{itemize}
      \item[1.] $x$ and $\Xi$ returned to initial values -- this point gives a trivial 
      transformation from the monodromy group. 
       \item[2.] The value of $x$ did not return to the initial values (within a given precision)
       -- continue integrating around this loop. If $x$ does not return to the initial values 
       after a sufficiently large number of loops, analyze the density of the trajectory 
       in the phase space (related to the Painlev\'e property). 
      \item[3.] The values of $x$ returned to the initial values after a finite number of loops, 
      but of $\Xi$ did not -- store the matrix $\Xi$, it is one of the generators of the monodromy group.
     \end{itemize}
  \item[iv.] Compute the pairwise commutators of all the matrices obtained in iii.3. 
  If there are non-vanishing commutators make a conclusion about non-integrability; if not
  choose another initial value of $x(t_0) = x_0$ in step iii. 
\end{enumerate}

\section{Application}

In this section we apply the developed algorithm to some systems
having mechanical interest. We start with the example that served one of the motivations to favor the approach of Morales--Ramis in comparison to the Ziglin's one -- the system described by the Henon-Heiles 
hamiltonian. Within the framework of our algorithm we can also study 
the systems we were interested in while describing another approach in \cite{int-num}, namely the 
pendulum-type systems and satellite dynamics. At the end we give some details of the implemented 
algorithm.

\subsection{Henon-Heiles system}
The Henon-Heiles hamiltonian describes a very simple model of 
a star moving close to the galactic center. It reads
$$
  H_h = \frac{1}{2}(p_1^2 + p_2^2) - q_2^2(A + q_1) - \frac{\lambda }{3} q_1^3.
$$
In \cite{audin} (referring also to original works \cite{morales,morales-ramis}) 
the monodromy group and the differential Galois group were constructed for the 
hamiltonian equations governed by $H_h$ using a rather simple particular solution 
satisfying $q_2=p_2 \equiv 0$. For this solution the monodromy group is trivial and does not 
obstruct integrability, while for $\lambda = 0$ the theory of Morales--Ramis 
permits to show that the system is non-integrable.
But for $\lambda \neq 0$ and $A \neq 0$ one can only state non-integrability 
for $\frac{6}{\lambda} \neq \frac{k(k+1)}{2}, \quad k \in \mathbb{Z}$.
For $A = 0$ the approach does not give any result. Later Morales and Ramis 
performed a more detailed study, showing that the question of 
integrability is open only for $\lambda = 1, 2, 6, 16$. 
Using the fact that one is not forced to be restricted to the above particular solution 
in our numerical approach, we can study the remaining cases by our method using a more 
complicated one. 

For example for $\lambda = 1, \quad A = 0.25$  consider the initial data
$(q_1, q_2, p_1, p_2) = (1, -0.4, -1.25, -0.3), \quad t_0 = 1$ 
(the numbers are chosen arbitrarily not to obtain a vanishing solution for $q_2, p_2$).
The commutator of matrices obtained by going along a loop around the points  
$(0.2 + 2.5i)$ and $(0.2 - 2.5i)$ is equal to 

$$ 
 \left(
\begin{array}{cccc}
0.39 - 0.91i & 0.66 - 3.46i & -1.10 + 2.57i & 0.66 - 2.22i\\
-0.94 + 2.30i & -0.64 + 2.60i &1.31 - 1.87i &-1.72 + 7.99i\\
0.31 - 0.76i &0.27 - 1.21i & -0.52 + 0.88i & 0.58 - 2.51i  \\
0.46 - 1.08i &0.80 - 4.20i &-1.36 + 3.13i &0.77 - 2.56i  
\end{array}
 \right),
$$
that results in non-integrability of the system.

\subsection{Triple pendulum}
A triple pendulum is a system of three mass-points (described by the radius-vectors
${\mathbf r}_i$) connected by weightless inextensible rods. 
We consider a free planar motion of this system which within the Lagrangian formalism 
can be described by a system with constraints of the form
\begin{eqnarray}
{\mathbf r}_1^2 - l_1^2 = 0,  \label{constr1} \nonumber
\\
({\mathbf r}_i - {\mathbf r}_{i-1})^2 - l_i^2 = 0,\qquad i = 2,3.
\nonumber
\end{eqnarray}
We are not going to describe the formalism in full details here since we have 
already sketched it together with the motivations in \cite{int-num} (cf. also references therein). Let us only recall 
that using a convenient parametrization by angles $\beta_1, \beta_2$ between the segments of 
the pendulum the system can be reduced by Routh transform to two degrees of freedom.
One technical difficulty is that we need to have an explicit form 
of the system of variational equations, which is in this case much more 
complicated than in the previous example. Luckily, we can use algorithms 
of symbolic computation to obtain it; in this case we used the Sage (\cite{sage}) software 
package.

Turning to the results, for the initial data
$\, (\beta_1, \dot \beta_1, \beta_2, \dot \beta_2) = 
(0.3, -1, -0.15 ,0.5 ),$ $t_0 = 1$, 
the commutator of the matrices obtained by going along the loop around the points 
$(0.2 + 0.5i)$ and $(0.2 - 0.5i)$ (each of them six times),
reads
$$ 
 \left(
\begin{array}{cccc}
0.62 - 0.62i & 0.81 + 0.83i & -0.24 - 0.34i & -0.02 + 0.33i\\ 
-0.07 - 0.14i & 1.04 + 0.10i & -0.04 + 0.20i & -0.02 - 0.39i \\ 
0.02 - 1.01i & 0.25 + 1.52i & 0.15 + 0.57i & 0.86 - 1.44i \\ 
-0.01 + 0.02i & -0.02 - 0.07i & -0.02 - 0.09i & 1.03 + 0.11i
\end{array}
 \right).
$$
That is the system is meromorphically non-integrable which is in perfect agreement with the 
results of \cite{int-num}.

\subsection{Satellite dynamics}
Let us now consider another example already mentioned in \cite{int-num} -- the 
motion of a dynamically symmetric satellite along a circular orbit (\cite{bardin}).
Using again the Routh reduction procedure one can describe the dynamics by the hamiltonian
$$
  H = \frac{\displaystyle p_{\psi}^2}{\displaystyle 2 \sin^2\theta} + \frac{\displaystyle p_{\theta}^2}{\displaystyle 2}
   - p_{\psi} + \frac{1}{2} \sin^2\psi \sin^2 \theta.
$$
Consider the trajectory starting at $\psi = 0, \theta = 1, p_{\psi} = 0.1, p_{\theta} = 0, t_0 = 0$
and continue it along the loop around the points $t = 4.8 + 0.8i$ and $t =4.8 - 0.8i$ 
(two times around each of them).
The respective commutator reads
$$
\left(
\begin{array}{llll}
8849.8 + 13.3i & 37.9 - 126.1i & 2044.5 + 35.4i & -1843.3 - 125.9i \\
-9456.3- 239.5i &311.8 - 62.7i &-2350.7 - 37.1i &1972.3 + 197.3i \\
-34540.9 - 527.5i &596.4 - 53.2i &-8340.9 - 82.1i &7205.5 + 624.6i \\
4032.6 - 556.4i &1615.1 - 971.5i &177.9 + 135.1i &-820.7 + 131.5i \\
\end{array}
\right).
$$
That is the system is also meromorphically non-integrable.
An interesting observation here is that the configuration of the system is parametrized 
periodically by the angles $\psi, \theta$, that is we didn't have to make them return exactly 
to the initial condition, but only modulo $2\pi$.
Another interesting feature of the system revealed by the numerical 
experiment, is that the complexified dynamics of it is much more sophisticated than the 
real one. That is if the system is indeed locally integrable as we have conjectured in 
\cite{int-num}, these integrals can not be continued to meromorphic functions 
in the complex domain.

\subsection{Implementation details}
\begin{center}
 \begin{figure}[htp] \centering
     \includegraphics*[width=0.9\linewidth]{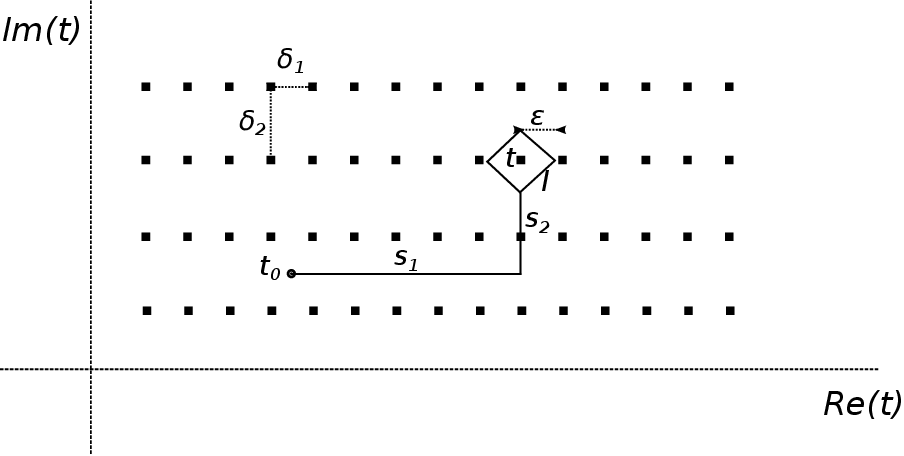}
    \caption{              \label{fig:map}
		 Grid of points $t$ and loops on the plane of complex time.
            }
 \end{figure}
\end{center}

Let us be more precise about some details that are related to fixing the 
freedom in the presented algorithm. 
For applying it one needs to choose the strategy for searching the 
branching points (steps ii. and iii.). In the above examples the loops 
around any point $t$ were always taken to be of the form shown on the schematic
figure \ref{fig:map}, 
i.e. the path starts from the distinguished point $t_0$, goes parallel first to the 
real axis (segment $s_1$), then to the imaginary one
(segment $s_2$), until it comes $\varepsilon$ close to 
the point $t$, after making a small loop around it ($l$) it goes back 
to $t_0$ again parallel to the axes. The points $t$ were chosen 
from a rectangular lattice with the parameters $\delta_1, \delta_2$. It is sufficient 
to find a couple of them that exhibit the behaviour of solution corresponding to 
the step iii.3. of the algorithm.
Then by adjusting $\varepsilon$ and the timestep one obtains the desired 
accuracy of numerical integration. 
In the above examples this permitted to obtain the Runge estimate 
of the error of numerical integration beyond the written digits, i.e. 
beyond the second digit after the decimal point (first one in the satellite example
due to high gradient norm). This permits to conclude that the commutators 
are not vanishing. Let us also note that in the presented examples the 
problem of resonance of the monodromy transformations is easy to handle:
it is sufficient to verify that the modulus of two eigenvalues of the computed monodromy 
matrices is different from $1$, which is again an open condition. 

It is sometimes useful to visualize the solution of (\ref{2syst1})
around the initial point $(t_0, x_0)$ by plotting the norm of the gradient 
of its flow along the paths containing just the segments $s_1$ and 
$s_2$ with a small parameter $\delta_1$ of the lattice: 
the numerical experiment shows that the maxima of this norm 
often correspond to branching points.

\section{Conclusion}
Thus, we have presented an algorithm for analysis of integrability 
of dynamical systems via the Ziglin's method using basically the 
properties of the monodromy group. 
As the examples show, it permits to extend the range of applicability of the method
mainly because it resolves the problem of finding an explicit particular solution
of the system. An important feature of the algorithm is that the trajectory obtained numerically, 
which is finally used for conclusion, is rather short and therefore can be computed with 
any given precision. It means that the algorithm indeed  provides a rigorous method 
of computer assisted proof of non-integrability. 

We have considered some mechanical examples that were not studied before, but 
the algorithm also has purely mathematical value. Namely, it is well adapted to 
computation of the monodromy group (or at least some subgroup of it) 
not necessarily related to integrability problem.
We expect this to be useful also in analysis of the complexity of the differential Galois group
which is the Zariski closure of the monodromy group, since 
by applying the algorithm one is able to produce the ``lower bound'' 
for it. 

To conclude, let us mention that one of the motivations for developing constructive
numerical methods of studying integrability for us is the qualitative analysis of the 
dynamical systems with delay or self-control appearing naturally in celestial 
mechanics and biological modeling. 

\begin{acknowledgments}
This work has been mostly done in the Claude Bernard Lyon 1 University and 
the Dorodnitsyn Computing Center of Russian academy of sciences -- the author is thankful to Sergey Stepanov for constant attention. The author also thanks Alexei Tsygvintsev for useful discussions at early 
stages of this work and the anonymous referee for important remarks clarifying the presentation.
\end{acknowledgments}



\end{document}